\documentclass[12pt]{amsart}
\usepackage{amsmath,amsthm,epsfig,amssymb}
\usepackage{amsfonts,amsmath}
\title{Primes of the form $\pm a^2\pm qb^2$}
\author{Eugen J. Ionascu and Jeff Patterson}
\curraddr{ Department of Mathematics\\ Columbus State University\\4225 University Avenue\\
Columbus, GA 31907} \email{math@ejonascu.edu, j3phr3y@gmail.com }   \subjclass{}
\date{April $7^{th}$, 2013}
\begin{document}
\def\sms{\small\scshape}
\baselineskip18pt
\newtheorem{theorem}{\hspace{\parindent}
T{\scriptsize HEOREM}}[section]
\newtheorem{proposition}[theorem]
{\hspace{\parindent }P{\scriptsize ROPOSITION}}
\newtheorem{corollary}[theorem]
{\hspace{\parindent }C{\scriptsize OROLLARY}}
\newtheorem{lemma}[theorem]
{\hspace{\parindent }L{\scriptsize EMMA}}
\newtheorem{definition}[theorem]
{\hspace{\parindent }D{\scriptsize EFINITION}}
\newtheorem{problem}[theorem]
{\hspace{\parindent }P{\scriptsize ROBLEM}}
\newtheorem{conjecture}[theorem]
{\hspace{\parindent }C{\scriptsize ONJECTURE}}
\newtheorem{example}[theorem]
{\hspace{\parindent }E{\scriptsize XAMPLE}}
\newtheorem{remark}[theorem]
{\hspace{\parindent }R{\scriptsize EMARK}}
\renewcommand{\thetheorem}{\arabic{section}.\arabic{theorem}}
\renewcommand{\theenumi}{(\roman{enumi})}
\renewcommand{\labelenumi}{\theenumi}

\newcommand{\comments}[1]{}
\def\phi{\varphi}
\def\ra{\rightarrow}
\def\sd{\bigtriangledown}
\def\ac{\mathaccent94}
\def\wi{\sim}
\def\wt{\widetilde}
\def\bb#1{{\Bbb#1}}
\def\bs{\backslash}
\def\cal{\mathcal}
\def\ca#1{{\cal#1}}
\def\Bbb#1{\bf#1}
\def\blacksquare{{\ \vrule height7pt width7pt depth0pt}}
\def\bsq{\blacksquare}
\def\proof{\hspace{\parindent}{P{\scriptsize ROOF}}}
\def\pofthe{P{\scriptsize ROOF OF}
T{\scriptsize HEOREM}\  }
\def\pofle{\hspace{\parindent}P{\scriptsize ROOF OF}
L{\scriptsize EMMA}\  }
\def\pofcor{\hspace{\parindent}P{\scriptsize ROOF OF}
C{\scriptsize ROLLARY}\  }
\def\pofpro{\hspace{\parindent}P{\scriptsize ROOF OF}
P{\scriptsize ROPOSITION}\  }
\def\n{\noindent}
\def\wh{\widehat}
\def\eproof{$\hfill\bsq$\par}
\def\ds{\displaystyle}
\def\du{\overset{\text {\bf .}}{\cup}}
\def\Du{\overset{\text {\bf .}}{\bigcup}}
\def\b{$\blacklozenge$}

\def\eqtr{{\cal E}{\cal T}(\Z) }
\def\eproofi{\bsq}

\begin{abstract}   Representations of primes by simple quadratic forms, such as $\pm a^2\pm qb^2$, is a subject that goes back to Fermat,
Lagrange, Legendre, Euler, Gauss and many others. We are interested in a comprehensive list of such results, for $q\le 20$.
Some of the results can be established with elementary methods and we put them at work on some instances.
We are introducing new relationships between various representations.
\end{abstract} \maketitle
\section{INTRODUCTION}

Let us consider the following three types representations of natural number:

\begin{equation}\label{firstrepr} {\cal E}(q):=\{n\in \mathbb N |n=a^2+qb^2, \text{with}\ a,b\in \mathbb Z\},
\end{equation}

\begin{equation}\label{secondrepr} {\cal H_1}(q):=\{n\in \mathbb N |n=qb^2-a^2, \text{with}\ a,b\in \mathbb Z\}, \ \text{and}
\end{equation}

\begin{equation}\label{thirdrepr} {\cal H_2}(q):=\{n\in \mathbb N |n=a^2-qb^2, \text{with}\ a,b\in \mathbb Z\}.
\end{equation}
\n We are going to denote by $\cal P$ the set of prime numbers.
In this paper we are going to exemplify how standard elementary methods can be used to obtain the representations stated in the next three
theorems:

\begin{theorem}\label{fermatstyle} For a prime $p$ we have $p\in {\cal E}(q)$ if and only if

{\tiny ({\bf I})} (Fermat) ($q=1$) $p=2$ or  $p\equiv 1$  (mod 4)

{\tiny ({\bf II})} (Fermat)  ($q=2$)   $p=2$ or  $p\equiv 1$ or $p\equiv 3$  (mod 8)

{\tiny ({\bf III})} (Fermat-Euler) ($q=3$) $p=3$ or $p\equiv 1$  (mod 6)

{\tiny ({\bf IV})} ($q=4$) $p\equiv 1$  (mod 4)

{\tiny ({\bf V})} (Lagrange) ($q=5$)  $p=5$ or $p\equiv j^2$  (mod 20) for some $j\in \{1, 3\}$

{\tiny ({\bf VI})}  ($q=6$)  $p\equiv 1$ or $7$  (mod 24)

{\tiny ({\bf VII})} ($q=7$)  $p=7$ or $p\equiv j^2$  (mod 14) for some $j\in \{1, 3, 5\}$

{\tiny ({\bf VIII})}  ($q=8$)  $p\equiv 1$  (mod 8)

{\tiny ({\bf IX})} ($q=9$)  $p\equiv j^2$  (mod 36) for some $j\in\{1,5,7\}$

{\tiny ({\bf X})} ($q=10$)  $p\equiv j$ (mod 40) for some $j\in\{1,9,11,19\}$

{\tiny ({\bf XI})}  ($q=12$) $p\equiv j$  (mod 48) for some $j\in\{1, 13, 25, 37\}$

{\tiny ({\bf XII})} ($q=13$) $p\equiv j^2$ (mod 52) for some $j\in\{1,3,5,7,9,11\}$

{\tiny  ({\bf XIII})}  ($q=15$)  $p\equiv j$  (mod 60) for some $j\in\{1, 19, 31, 49\}$

{\tiny  ({\bf XIV})}  ($q=16$)  $p\equiv 1$  (mod 8)

\end{theorem}

We are going to prove {\small (VII)}, in order introduce the method that it will be employed several times.
One may wonder what is the corresponding characterization for $q=11$ or $q=14$. It turns out that an answer cannot be
formulated only in terms of residue classes as shown in  (\cite{sw1992}). We give in Theorem~\ref{newcharcterzations} possible characterizations
whose proofs are based on non-elementary techniques which are described in  \cite{dcox}.

\begin{theorem}\label{dualpell} For a  prime $p$ we have $p\in {\cal H_1}(q)$ if and only if

{\tiny ({\bf I})} ($q=1$) $p\not =2 $

{\tiny ({\bf II})} ($q=2$) $p=2$ or $p\equiv \pm 1$  (mod 8) (i.e. $p\equiv 1$  or  $p\equiv 1$ (mod 8))

{\tiny ({\bf III})} ($q=3$) $p\in\{2,3\}$ or $p\equiv 11$  (mod 12)

{\tiny ({\bf IV})} ($q=4$) $p\equiv 3$  (mod 4)

{\tiny ({\bf V})} ($q=5$)  $p=5$ or $p\equiv \pm j^2$ (mod 20) for some $j\in\{1,3\}$

{\tiny ({\bf VI})}  ($q=6$)  $p=2$ or $p\equiv j$ (mod 24) for some $j\in\{5,23\}$

{\tiny ({\bf VII})} ($q=7$)  $p=7$ or $p\equiv j$  (mod 14) for some $j\in \{3, 5, 13\}$

{\tiny ({\bf VIII})}  ($q=8$)  $p=7$ or $p\equiv -j^2$  (mod 32) for some $j\in \{1,3, 5,7\}$

{\tiny ({\bf IX})} ($q=9$)  $p\equiv -1$  (mod 6)

{\tiny ({\bf X})} ($q=10$)  $p\equiv j$ (mod 40) for some $j\in\{1,9,31,39\}$

{\tiny ({\bf XI})}  ($q=11$) $p\in \{2,11\}$ or $p\equiv -j^2$  (mod 44) for some $j\in\{1, 3, 5, 7,9\}$
\end{theorem}

In this case, for exemplification, we show {\small (V)}.

\begin{theorem}\label{pell} For a  prime $p$ we have $p\in {\cal H_2}(q)$  if and only if

{\tiny ({\bf I})} ($q=1$) $p\not =2 $

{\tiny ({\bf II})} ($q=2$) $p=2$ or $p\equiv \pm 1$  (mod 8)

{\tiny ({\bf III})} ($q=3$)   $p\equiv 1$  (mod 12)

{\tiny ({\bf IV})} ($q=4$) $p\equiv 1$  (mod 4)

{\tiny ({\bf V})} ($q=5$) $p=5$ or  $p\equiv \pm j^2$ (mod 20) for some $j\in\{1,3\}$

{\tiny ({\bf VI})}  ($q=6$)  $p=3$ or $p\equiv j$ (mod 24) for some $j\in\{1,19\}$

{\tiny ({\bf VII})} ($q=7$)  $p=2$ or $p\equiv j$  (mod 14) for some $j\in \{1, 9, 11\}$

{\tiny ({\bf VIII})}  ($q=8$)  $p=7$ or $p\equiv j^2$  (mod 32) for some $j\in \{1,3, 5,7\}$

{\tiny ({\bf IX})} ($q=9$)  $p\equiv 1$  (mod 6)

{\tiny ({\bf X})} ($q=10$)  $p\equiv j$ (mod 40) for some $j\in\{1,9,31,39\}$

{\tiny ({\bf XI})}  ($q=11$) $p\equiv j^2$  (mod 44) for some $j\in\{1, 3, 5, 7,9\}$
\end{theorem}

We observe that for $q=2$, $q=5$, $q=10$ the same primes appear for both characterizations in Theorem~\ref{dualpell} and Theorem~\ref{pell}.
There are several questions that can be raised in relation to this observation:

{\bf Problem 1:} Determine all values of $q$, for which we have

\begin{equation}\label{samerepresentation}
{\cal H_1}(q)\cap {\cal P}={\cal H_2}(q)\cap {\cal P}.
\end{equation}

{\bf Problem 2:} If the equality (\ref{samerepresentation}) holds true for relatively prime numbers $q_1$ and $q_2$, does is it hold true for $q_1q_2$?

In \cite{dcox}, David Cox begins his classical book  on the study of (\ref{firstrepr}), with a detailed and well documented historical introduction of the main ideas used and the difficulties encountered in the search of new representations along time.  The following abstract characterization in \cite{dcox}   brings more light into this subject:

\begin{quote}
    {\bf (Theorem 12.23 in \cite{dcox})} {\em Given a positive integer $q$, there exists an irreducible polynomial with integer coefficients $f_q$ of degree $h(-4q)$, such that  for every odd  prime $p$ not dividing $q$,
    $$p=a^2+qb^2 \Leftrightarrow \text{the equations} \begin{cases} x^2\equiv -q \ (mod\ p) \\ \\ f_q(x)\equiv 0\ (mod\ p)
    \end{cases}$$ have integer solutions. An algorithm for computing $f_q$ exists. ($h(D)$ is the number of classes of primitive positive definite quadratic forms of discriminant $D$).}
\end{quote}

\n     While some of the representations included here are classical and other may be more or less known. We found some of the polynomials included here
by computational experimentations. For more details in this direction see \cite{bc}, \cite{db},  \cite{hc2}, \cite{dcox}, \cite{hw1}, \cite{fl} and \cite{sw1992}.

\begin{theorem}\label{newcharcterzationspolynomials} For an odd prime $p$ we have $p=a^2+qb^2$ for some integers $a$, $b$ if and only if

{\tiny ({\bf I})}  ($q=11$)  $p>2$ and the equation
 $$(X^3+2X)^2+44 \equiv 0\ (mod \ p)\ \text{ has a solution,}$$

{\tiny  ({\bf II})}  (Euler's conjecture) ($q=14$)  the equations
  $$X^2+14\equiv 0\   and \ (X^2+1)^2-8\equiv 0\ (mod \ p) \ \text{have solutions}$$

({\tiny  {\bf III}}) ($q=17$) the equations $X^2+17\equiv 0$ and $(X^2-1)^2+16\equiv 0$ (mod p) have solutions

{\tiny  ({\bf IV})}  ($q=18$) the equation $(X^2-3)^2+18(2^2)\equiv 0$ (mod p) has a solution

{\tiny  ({\bf V})} ($q=19$) the equation $(X^3-4x)^2+19(4^2)\equiv 0$ (mod p) has a solution

{\tiny ({\bf VI})} ($q=20$)   the equation  $(X^4-4)^2+20X^4\equiv 0$ (mod p)  has a solution

{\tiny  ({\bf XXI})}  ($q=21$)  the equation  $(X^4+4)^2+84X^4\equiv 0$ (mod p)  has a solution

{\tiny  ({\bf XXII})}  ($q=22$) $p>22$ and the equation  $(x^2+3)^2+22(4^2)\equiv 0 $ (mod p)  has a solution

{\tiny   ({\bf XXIII})}  ($q=23$)  the equation $(X^3+15X)^2+23(19^2)\equiv 0$ (mod p) has a solution

{\tiny   ({\bf XXIV })}  ($q=24$) the equation $(X^4+4)^2+24(2X)^4\equiv 0$ (mod p) has a solution

{\tiny    ({\bf XXV }) }  ($q=25$) $p>25$ the equation $X^4+100 \equiv 0$ (mod p) has a solution

{\tiny     ({\bf XXVI })}  ($q=26$)

 {\tiny ({\bf XXVII})} (Gauss)  ($q=27$)  $p\equiv 1$    (mod 3) and the equation $X^3\equiv 2$ (mod $p$) has a solution;

{\tiny  ({\bf XXVIII})}  ($q=28$)

{\tiny  ({\bf XXVIV})}  ($q=29$) $p\equiv 1$    (mod 4) and  the equation $(X^3-X)^2+116=0$ (mod $p$)  has a solution;

 {\tiny  ({\bf XXX})}  ($q=30$)

{\tiny  ({\bf XXXI})} ($q=31$)  (L. Kronecker, pp. 88 \cite{dcox}) the equation
$$(X^3-10X)^2+31(X^2-1)^2\equiv 0\ (mod \ p)\ \text{ has a solution}$$

{\tiny  ({\bf XXXII})} ($q=32$)  $p\equiv 1$    (mod 8) and the equation

$$(X^2-1)^2\equiv -1\ (mod \ p) \ \text{has a solution.}$$

{\tiny  ({\bf XXXVII})} ($q=37$) the equation  $X^4+31X^2+9=0$ (mod $p$) has a solution

{\tiny ({\bf LXIV})} (Euler's conjecture)  ($q=64$)  $p\equiv 1$    (mod 4) and the equation $X^4\equiv 2$ (mod $p$) has a solution.
\end{theorem}

 Our interest in this subject came from studding the problem of finding all equilateral triangles, in the three dimensional space, having integer coordinates for their vertices (see \cite{rceji}, \cite{eji}, \cite{ejirt}, and \cite{ejips}). It turns out that such equilateral triangles exist only in planes  ${\cal P_{a,b,c,f}}:=\{(x,y,z)\in \mathbb R^3:ax+by+cz=f,f\in \mathbb Z\}$ where $a$, $b$, and $c$ are in such way

\begin{equation}\label{first}
a^2+b^2+c^2=3d^2
\end{equation}

\n for some integer $d$ and side-lengths of the triangles are of the form $$\ell=d\sqrt{2(m^2-mn+n^2)}$$
\n for some integers $m$ and $n$. This leads to investigations of primes of the first three forms in the  Let us include here a curious fact that we ran into at that time (\cite{eji}).

\begin{proposition}\label{thebeauty}  An integer $t$ which can be written as  $t=3x^2-y^2$ with
$x,y \in \mathbb Z$ is the sum of two squares if and only if $t$ is of
the form $t=2(m^2-mn+n^2)$ for some integers $m$ and $n$.
\end{proposition}
\n If we introduce the sets $A:=\{t\in \mathbb Z| t=3x^2-y^2, x,y \in \mathbb Z\}$, $B:=\{t\in \mathbb Z| t=x^2+y^2, x,y \in \mathbb Z\}$ and
$C:=\{t\in \mathbb Z| t=2(x^2-xy+y^2), x,y \in \mathbb Z\}$ then we actually have an interesting relationship between these sets.

\begin{theorem}\label{trinity} For the sets defined above, one has the inclusions
\begin{equation}\label{inclusions}A\cap B\subsetneqq C,\ \  B\cap C \varsubsetneqq A, \ \ \text{and}\  \ C\cap A\varsubsetneqq B.
\end{equation}
\end{theorem}

\begin{figure}
\[
\underset{  }{\epsfig{file=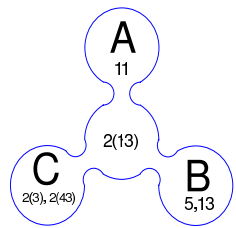,height=1in,width=1in}}
\]
\caption{\em ``God created the integers, all else is the work of man." Leopold Kronecker} \label{GOD}
\end{figure}

\n  We include a proof of this theorem in the Section~\ref{trinitysection}. The inclusions in (\ref{inclusions}) are strict as one can see from Figure~\ref{GOD}.

Let us observe that there are primes $p$ with the property that $2p$ is in all three sets $A$, $B$ and $C$.
We will show that these primes are the primes of the form $12k+1$ for some integer $k$. Some  representations for such primes are included next:
\begin{equation}\label{trinityprimes}\begin{array}{c}  13=(1^2+5^2)/2=3^2-3(4)+4^2=[3(3^2)-1]/2\\ \\
37=(5^2+7^2)/2=3^2-3(7)+7^2=[3(5)^2-1]/2\\ \\
61=(1^2+11^2)/2=4^2-4(9)+9^2=[3(9^2)-11^2]/2.
\end{array}
\end{equation}

It is natural to ask whether or not  the next forms in the Theorem~\ref{fermatstyle} aren't related to similar parameterizations for regular or semi-regular simplices in $\mathbb Z^n$ for bigger values of $n$. In \cite{s}, Isaac Schoenberg gives a characterization of those $n$'s for which a regular simplex exists in $\mathbb Z^n$. Let us give the  restatement of Schoenberg's result which appeared in \cite{mc}: {\it all $n$ such that $n+1$ is a sum of 1, 2, 4 or 8 odd squares.}

We summarize next the explicit forms of Theorem 12.23 in \cite{dcox} for $q\le 100$.

As interesting corollaries of these statements we see that if one prime $p$ has some representation it must have some other type of representation(s).
Let us introduce a notation for these classes of primes: $$\cal P_q:=\{p \ \text{ odd prime}|p=a^2+qb^2 \ \text{for some} \ a,b\in \mathbb N\}.$$
So we have $\cal P_1=\cal P_{4}$,  $\cal P_8=\cal P_{16}$ (Gauss, see \cite{jvu}), $\cal P_5\subset \cal P_1$, $\cal P_{10}\subset \cal P_{2}$, .....
In the same spirit, we must bring to reader's attention, that in the case $q=32$ there exists a characterization due to Barrucand and Cohn \cite{bc}, which can be written with our notation as

$$\cal P_{32}=\{p\ |\ p\equiv 1\ (mod\ 8), \text{there exists}\ x\ \text{such that}\ x^8\equiv -4 \ (mod \ p)\}.$$

We observe that ({\bf xx}) implies this characterization because $x^8+4=(x^4-2x^2+2)(x^4+2x^2+2)$ and clearly $(x^2-1)^2+1=x^4-2x^2+2$. In fact, the two statements are equivalent. Indeed,
 if $a$ is a solution of $x^8+4\equiv 0$ (mod p) then we either have $x^4-2x^2+2\equiv 0$ (mod p ) or $x^4+2x^2+2\equiv 0$ (mod p). We know that there exists a solution $b$  of
 $x^2+1\equiv 0$ (mod p). Hence if $a^4+2a^2+2\equiv 0$ (mod p) then $(ab)^4-2(ab)^2+2\equiv 0$ (mod p) which shows that the equation $x^4-2x^2+2\equiv 0$ (mod p ) always has a solution.

Also, another classical result along these lines is Kaplansky's Theorem (\cite{ik}):

\begin{theorem}\label{kaplansky} A prime of the form $16n+9$ is in $\cal P_{32}\setminus \cal P_{64}$ or in $\cal P_{64}\setminus \cal P_{32}$.
For a prime $p$ of the form $16n+1$ we have  $p\in \cal P_{32}\cap \cal P_{64}$ or $p\not \in \cal P_{64}\cup \cal P_{32}$.
\end{theorem}

For further developments similar to Kaplansky's result we refer to \cite{db}.
One can show that the representations in Theorem~\ref{fermatstyle} are unique (see Problem 3.23 in \cite{dcox}).

\section{Case {\bf (vii)}}

We are going to use elementary ideas in the next three sections.

\begin{theorem}{\bf [Gauss]}\label{lawofquadraticreciprocity}
For every $p$ and $q$ odd prime numbers we have

\begin{equation}\label{lqr}
\left(\frac{p}{q}\right)\left(\frac{q}{p}\right)=(-1)^{\frac{p-1}{2}\frac{q-1}{2}}.
\end{equation}

\n with notation $\left(\frac{\cdot}{p}\right)$, defined for every odd
prime $p$ and every $a$ coprime with $p$ known as the Legendre
symbol:

\begin{equation}\label{legendresymbol}
\left(\frac{a}{p}\right)=\begin{cases} 1 \ if\ the \ equation\
x^2\equiv a\ (mod\ p)\ has\ a\ solution, \\ \\ -1\ if\ the \
equation\ x^2\equiv a\ (mod\ p)\  has\ no \ solution
\end{cases}
\end{equation}
\end{theorem}
We think that this method can be used to prove all the statements in Theorem~\ref{fermatstyle} except (xi), (xiv), and  (xvii)-(xx). The exceptions are either known (see \cite{dcox}) or shown in the last section.  We learned about this next technique from \cite{nzm1991} and \cite{r2005}.

\n {\it Necessity:} If $p=x^2+7y^2$ then  $p\equiv x^2$ (mod 7). Clearly we may assume $p>7$. Therefore, $x$  may be assumed to be different of zero.
Then the residues of $p$ (mod 7) are 1, 2, or 4. Let us suppose that $p\equiv r$ (mod 14) with $r\in \{0,1,2,..., 13\}$.  Because $p$ is prime, $r$ must be an odd number, not a multiple of 7 and which equals 1, 2 or 4 (mod 7). This leads to only three such residues, i.e. $r\in \{1,9,11\}$, which are covered by the odd squares $j^2$, $j\in \{1,3,5\}$.

\n {\it Sufficiency:} We may assume that $p>2$. Let us use the hypothesis to show that the equation $x^2=-7$ has a solution.
Let $p$ be a prime of the form $14k+r$, $r\in \{1,9,11\}$, $k\in \mathbb N\cup \{0\}$.
By the Quadratic Reciprocity, we have $(\frac{7}{p})(\frac{p}{7})=(-1)^{\frac{3(p-1)}{2}}$.
Since $(\frac{-1}{p})=(-1)^{\frac{p-1}{2}}$, then $$(\frac{-7}{p})=(-1)^{\frac{p-1}{2}}(\frac{7}{p})=(-1)^{\frac{p-1}{2}+\frac{3(p-1)}{2}}(\frac{p}{7})=(\frac{r'}{7}), \ {\text where} \ p=7(2k')+r', r'\in \{0,1,..,6\}.$$
This shows that if $r'\in \{1,2,4\}$ we have a solution $x_0$ for the equation $x^2\equiv -7$ (mod p).

Let us now apply the Pigeonhole Principle: we let $m\in \mathbb N$ be in such a way that $m^2<p<(m+1)^2$.
We consider the function $g:\{0,1,2,...,m\}\times \{0,1,2,...,m\}\to \{0,1,2,....,p-1\}$ defined by $g(u,v)\equiv u+vx_0$ (mod p).
Since $(m+1)^2>p$, we must have two distinct pairs $(a'',b'')$ and $(a',b')$ such that $g(a'',b'')=g(a',b')$. Then $a''-a'\equiv (b'-b'')x_0 $ (mod p).
Then, if we let $a=a''-a'$, and $b=b'-b''$ we get that $0<q:=a^2+7b^2\equiv b^2(x_0^2+7)\equiv 0$ (mod p).
But, $q=a^2+7b^2\le m^2+7m^2=8m^2<8p$. It follows that $q\in \{p,2p,3p,4p,5p,6p,7p\}$. We need to eliminate the cases $q\in \{2p,3p,4p,5p,6p,7p\}$.
If $q=7p$ then $7p=a^2+7b^2$ which implies that $a$ is a multiple of $7$, or $a=7a'$, which gives $p=b^2+7a'^2$ as wanted.

If $q=3p$, then $q=3(14k'+r)=7\ell+s$ where $s\in \{3,5,6\}$. But this is impossible because $q\equiv a^2$ (mod 7).
The same argument works if $q=6p$, because $r'\in \{1,2,4\}$ if and only if $6r' \in \{3,5,6\}$ (mod 7). Similarly,
the case $p=5p$ is no difference.

If $q=2p$ or $a^2+7b^2=2p$ implies that $a$ and $b$ cannot be both odd, since in this case $a^2+7b^2$ is a multiple of $8$ and $2p$ is not.
Therefore $a$ and $b$ must be both even, but that shows that $2p$ is a multiple of $4$. Again this is not the case.

Finally, if $q=4p$ then the argument above works the same way but in the end we just simplify by a 4.\eproof

\section{Cases $q\in \{11, 17, 19\}$}

Let us observe that the characterizations in Theorem~\ref{newcharcterzationspolynomials} cannot be easily checked for big primes $p$. 
Next, we use still similar elementary methods to show the following fact which seems to be the best what one can hope for in terms of a characterization in which certain quadratic forms
of the form $a^2+qb^2$ cannot be separated by simply the quadratic residues of odd numbers modulo $4q$.

\begin{theorem}\label{newthingy} (i) A prime $p>17$ is of the form $a^2+17b^2$ or $2p=a^2+17b^2$, for some
$a,b\in \mathbb N$ if and only if $p\equiv (2j+1)^2$ (mod 68) for some $j=0,...,7$.\par
(ii) The representation of a prime as in part (a) is exclusive, i.e. a prime $p$ cannot be of the form $a^2+17b^2$ and at the same time $2p=x^2+17y^2$, for some
$x,y\in \mathbb N$.
\end{theorem}

\proof \ (i) \n {\it $``\Rightarrow"$} If the prime $p$ can be written $p=a^2+17b^2$ then  $p\equiv a^2$ (mod 17) with $a$ not divisible by $17$.
 We observe that $a$ and $b$ cannot be both odd or both even. Then $p\equiv 1$ (mod 4). If $p=68k+r$ with $r\in \{0,1,2,...,67\}$ then $r\equiv 1$ (mod 4), not a multiple of $17$ and
a quadratic residue modulo 17, i.e. $r=17\ell +r'$ with $r'\in \{1, 2, 4, 8, 9, 13, 15, 16\}$. This gives
$r\in \{ 1, 9, 13, 21, 25, 33, 49, 53\}$. One can check that these residues are covered in a one-to-one way by the odd squares $j^2$, $j\in \{1,3,5,7,9,11,13,15\}$.

If $2p=a^2+17b^2$ then $2p\equiv a^2$ (mod 17) with $a$ not divisible by $17$. In this case $a$ and $b$ must be both odd and then $2p=a^2+17b^2\equiv 2$ (mod 8). This implies, as before, that $p\equiv 1$ (mod 4). If  $p=68k+r$ with $r\in \{0,1,2,...,67\}$ then $r\equiv 1$ (mod 4), not divisible by $17$ and
$2r$ is a quadratic residue modulo 17. Interestingly enough, we still have $r\in \{ 1, 9, 13, 21, 25, 33, 49, 53\}$.

\n {\it $``\Leftarrow"$} We have $p\equiv j^2$ (mod 17) and so $(\frac{p}{17})=1$. By the Theorem~\ref{lawofquadraticreciprocity}, we have $(\frac{17}{p})(\frac{p}{17})=(-1)^{8\frac{(p-1)}{2}}=1$ which implies $(\frac{17}{p})=1$.

Since $(\frac{-1}{p})=(-1)^{\frac{p-1}{2}}$, we get that $(\frac{-17}{p})=(-1)^{\frac{p-1}{2}}$. If $p=68k+j^2$ with $j\in \{1,3,5,7,9,11,13,15\}$,
we see that $(\frac{-17}{p})=1$. Therefore $x^2\equiv -17$ (mod p) has a solution $x_0$. As in the case $q=7$, if we use the same idea of the Pi
Pigeonhole Principle we obtain that $q=a^2+17b^2<18p$ for some $a,b\in \mathbb Z$ and $q\equiv 0$ (mod p). Hence $q=\ell p$ with $\ell \in \{1,2,...,17\}$.
We may assume that $gcd(a,b)=1$, otherwise we can simplify the equality $q=\ell p$ by  $gcd(a,b)$ which cannot be $p$.
Clearly if $\ell=1$, $\ell=2$ or $\ell=17$ we are done. Since $q\equiv 0$, 1 or 2 (mod 4) and $p\equiv 1$ (mod 4) we cannot have $\ell\in \{3,7,11,15\}$.
If $\ell \in \{4,8,12,16\}$, $\ell=4\ell'$, we can simplify the equality by a $4$ and reduce this case to $\ell' \in\{1,2,3,4\}$. Each one of these situations
leads to either the conclusion of our claim or it can be excluded as before or reduced again by a $4$.

(Case $\ell=5$ or $\ell=10$) Hence $q=\ell p=a^2+17b^2\equiv a^2+2b^2\equiv 0$ (mod 5). If $b$ is not a multiple of $5$ then this implies $x^2\equiv -2$ (mod 5) which is not true.
Hence $b$ must be a multiple of $5$ and then so must be $a$. Then the equality $\ell p=a^2+17b^2$ implies that $\ell p$ is a multiple of $25$ which is not possible.

(Case $\ell=6$ or $\ell=14$) In this case we must have $a$ and $b$ odd and then $q=2(4s+1)=\ell p$ which is not possible.

(Case $\ell=13$) In this case   $4q=(2a)^2+17(2b)^2=2p(3^2+17(1)^2)$. We will use Euler's argument (\cite{dcox}, Lemma 1.4, p. 10) here.
If we calculate $M=(2b)^2[3^2+17(1)^2]-4q=[3(2b)-2a][3(2b)+2a]$, we see that $2(13)$ divides $M$ and so it divides either $3(2b)-2a$ or $3(2b)+2a$. Without loss of generality we may assume that $2(13)$ divides $3(2b)-2a$. Hence, we can write $3(2b)-2a=2(13)d$ for some $d\in \mathbb Z$. Next, we calculate

$$2a+17d=3(2b)-2(13)d+17d=3(2b)-9d,$$

\n which implies that $2a+17d=3e$ for some $e\in \mathbb Z$. Also, from the above equality we get that $2b=e+3d$. Then

$$2p(26)=4q=(2a)^2+17(2b)^2=(3e-17d)^2+17(e+3d)^2=26(e^2+17d^2)\Rightarrow$$

$$2p=e^2+17d^2.$$

(Case $\ell=9$) We have   $4q=(2a)^2+17(2b)^2=2p(1^2+17(1)^2)$.
We calculate $M=(2b)^2[1^2+17(1)^2]-4q=(2b-2a)(2b+2a)$, we see that $2(9)$ divides $M$ and so it divides either $2b-2a$ or $2b+2a$.
We need to look into two possibilities now. First $2(9)$ divides one of the factors $2b-2a$ or $2b+2a$, or $2(3)$ divides each one of them. In the second situation
we can see that $3$ divides $4a=2b+2a-(2b-2a)$ and so $3$ must divide $b$ too. This last possibility is excluded by the assumption that $gcd(a,b)=1$.
Without loss of generality we may assume that $2(9)$ divides $2b-2a$. Hence, we can write $2b-2a=2(9)d$ for some $d\in \mathbb Z$. We set, $2a=e-17d$ and observe that
$2b=2a+18d=e-17d+18d=e+d$.
Then
$$2p(18)=4q=(2a)^2+17(2b)^2=(e-17d)^2+17(e+d)^2=18(e^2+17d^2)\Rightarrow$$

$$2p=e^2+17d^2.$$

\n (ii) To show this claim, we may use Euler's argument as above. \eproof

For primes $q$ which are multiples of four minus one, the patterns suggest that we have to change the modulo but also there are more trickier changes. Let us look at the cases $q=11$ and $q=19$. In case $q=11$, we have seen that the quadratic form
$a^2+11b^2$ in Theorem~\ref{fermatstyle}  can be separated by a polynomial from the other possible forms of representing primes which are quadratic  residues of odd numbers modulo 22.

\begin{theorem}\label{newthingy11} (i) A prime $p>11$ is of the form $a^2+11b^2$ or $3p=a^2+11b^2$, for some
$a,b\in \mathbb N$ if and only if $p\equiv (2j+1)^2$ (mod 22) for some $j=0,...,4$.\par
(ii) A prime $p>19$ satisfies $4p=a^2+19b^2$, for some
$a,b\in \mathbb N$ if and only if $p\equiv (2j+1)^2$ (mod 38) for some $j=0,...,8$.\par

(iii) The representations of a prime as in part (i) are exclusive, i.e. a prime $p$ cannot be in both representations.
\end{theorem}

We leave these proofs for the interested reader.
\section{Proof of Theorem~\ref{trinity}}\label{trinitysection}

Clearly the inclusions $A\cap B\subset C$ and  $C\cap A\subset B$  are covered by Proposition~\ref{thebeauty}.  To show   $B\cap C \subset A$ we will first prove it for $t=2p$ with
$p$ a prime. Since $2p=a^2+b^2$ we have $a^2\equiv -b^2$ (mod $p$). Because $p>2$, $a$ cannot be divisible by $p$ and so it has an inverse (mod $p$) say $a^{-1}$. This shows that $x_0=ba^{-1}$ is a solution of the equation $x^2\equiv -1$ (mod $p$). Similarly since $2p=2(x^2-xy+y^2)$ we get that $4(x^2-xy+y^2)=(2x-y)^2+3y^2\equiv 0$ (mod $p$).
This gives a solution $y_0$ of the equation  $x^2\equiv -3$ (mod $p$). So, we have $z_0=x_0y_0$ satisfying $z_0^2\equiv 3$ (mod $p$).
Let us now apply the Pigeonhole Principle as before: we let $m\in \mathbb N$ be in such a way that $m^2<p<(m+1)^2$.
We consider the function $g:\{0,1,2,...,m\}\times \{0,1,2,...,m\}\to \{0,1,2,....,p-1\}$ defined by $g(u,v)\equiv u+vz_0$ (mod $p$).
Since $(m+1)^2>p$, we must have two distinct pairs $(a'',b'')$ and $(a',b')$ such that $g(a'',b'')=g(a',b')$. Then $a''-a'\equiv (b'-b'')z_0 $ (mod $p$).
Then, if we let $r=a''-a'$, and $s=b'-b''$ we get that $q:=r^2-3s^2\equiv s^2(z_0^2-3)\equiv 0$ (mod $p$). So, $q$ needs to be a multiple of $p$. If $q=0$ then $r=\pm s\sqrt{3}$ which is not possible because $r$ and $s$ are integers not both equal to zero. If $q>0$ then $0<q\le r^2<p$, which is again impossible. It remains that $q<0$, and so $0<-q=3s^2-r^2\le 3s^2<3p$.
This leaves only two possibilities for $q$: either $q=-p$ or $q=-2p$. Hence, we need to exclude the case $3s^2-r^2=p$. This implies $4p=12s^2-4r^2=(2x-y)^2+3y^2$. Then $4r^2+(2x-y)^2\equiv 0$ (mod $3$). Since $-1$ is not a quadratic residue modulo $3$ we must have $r$ divisible by $3$ which is gives $p=3$ but we cannot have $6=a^2+b^2$. It remains that $2p=3s^2-r^2$.
Let us observe that in this case $s$ and $r$ cannot be both even or of different parities since $p$ must be of the form $4k+1$. Hence, we have the representation $p=(\frac{3s+r}{2})^2-3(\frac{s+r}{2})^2$.

To prove the inclusion in general we just need to observe that for any number $t\in B\cap C$ and a prime $p>2$ dividing $t$, then if $p$ is of the form $4k+3$ then it divides $a$ and $b$ and so $p^2$ divides $t$. The same is true if $p$ is of the form $6k-1$.  Clearly all the primes that appear in the decomposition of $t$ to an even power they can be factored out and reduce the problem to factors of the form $12k+1$ but for these factors we can apply the above argument and use the identities:

$$(y^2-3x^2)(v^2-3u^2)=(3ux+vy)^2-3(xv+uy)^2,$$

$$2(x^2-3y^2)=3(x+y)^2-(x+3y)^2.$$

\end{document}